\documentclass[12pt,reqno]{amsart}
\usepackage{graphicx} 
\def\CC{\mathbb C}

\def\PP{\mathbb P}

\def\RR{\mathbb R}
\def\NN{\mathbb N}
\def\Ibb#1{{\rm I\kern-.23em#1}}

\usepackage{eucal}
\renewcommand{\Re}{\operatorname{Re}}
\renewcommand{\Im}{\operatorname{Im}}
\renewcommand{\tilde}{\widetilde}

\DeclareMathOperator{\Res}{Res}
\DeclareMathOperator{\dist}{dist}
\numberwithin{equation}{section}

\newcommand{\bea}{\begin{eqnarray*}}
\newcommand{\eea}{\end{eqnarray*}}

\newtheorem{lem}{LEMMA}[section]
\newtheorem{theo}[lem]{THEOREM}

\newtheorem{coro}[lem]{COROLLARY}
\newtheorem{prop}[lem]{PROPOSITION}

\newtheorem{proposition}[lem]{PROPOSITION}

\begin{document}
\title[Bergman kernel on the intersection of two balls]
{The Bergman kernel on the intersection of two balls
in $\CC^2$}
\author{David E. Barrett and Sophia Vassiliadou}
\thanks{First author supported in part by the National Science Foundation
under Grant No. DMS-0072237. This work was completed while the second author 
was visiting the University of Toronto. She would like to thank the 
Department of Mathematics for its  hospitality and support.}
\address{Dept. of Mathematics\\University of Michigan
\\Ann Arbor, MI  48109-1109  USA}
\date{\today}
\email{barrett@umich.edu,\; sophia@math.lsa.umich.edu}
\begin{abstract}  We obtain an asymptotic expansion and 
some regularity results  for the 
Bergman kernel on the intersection of two balls in $\CC^2$.
 
\end{abstract}
\maketitle

%\newpage

 \section{Introduction}\label{S:intro}

Let $B_1=\{(z_1,z_2);\; |z_1|^2+|z_2|^2< 1\},\;\;B_2=\{(z_1,z_2);\;\;
|z_1-a_1|^2+|z_2-a_2|^2< r^2\}$ be two
balls in 
$\CC^2$ such that $\partial B_1, \;\partial B_2$ intersect real transversally.
There exist  two complex tangent points $p,q\in\partial B_1\cap \partial B_2$.

\medskip

In this paper we obtain an
asymptotic expansion for the Bergman kernel of  $\Omega:=B_1\cap B_2$.  (The Bergman kernel
function $K_\Omega$ on $\Omega\times\Omega$ is 
characterized by the conditions that
$\overline{K_\Omega(z,\cdot)}$ is holomorphic and square-integrable for all
$z$, and that $\int_\Omega K_\Omega(z,\zeta) f(\zeta)\,dV_\zeta=f(z)$ for all
holomorphic square-integrable $f$.  $K_\Omega$ satisfies $K_\Omega(z,\zeta) 
=\overline{K_\Omega(\zeta,z)}$.)

A generic description of our main result runs as follows.

\begin{theo}\label{T:maingen}

\begin{enumerate}
\item[i)]\label{I:tsym} For each $\zeta\in\Omega$ the Bergman kernel function
\linebreak
$K_\Omega(z,\zeta)$ is holomorphic in a neighborhood of $\overline
\Omega\setminus\{p,q\}$.
\item[ii)]  For $z$ near a complex tangent point $q$ we have an
asymptotic expansion of the form
\begin{equation}\label{E:crudeasym}
K_\Omega(z,\zeta)\sim
\sum_j \langle z-q, ^{T}a\rangle^{n_j}\langle z-q, a\rangle^{\gamma_j}
P_j\left(\log\langle z-q, a\rangle,\zeta\right):
\end{equation}
here $a=(a_1,a_2);$ $^{T}a=\left(\overline{a_2}, -\overline{a_1}\right);$ 
$<,>$ denotes the hermitian inner
product in 
$\CC^n;$ the
$n_j$ are nonnegative integers; the  $\gamma_j$ are (possibly complex)
exponents lying in the half-plane $\Re \gamma_j>-1-\frac{n_j}{2}$; and the
$P_j\left(\log\langle z-q, a\rangle,\zeta\right)$ are polynomials in $\log\langle
z-q, a\rangle$ with coefficients varying anti-holomorphically with $\zeta$.
\end{enumerate}
\end{theo}

See Theorem \ref{T:main} below for a more precise description of the expansion; up
to a change of coordinates it is equivalent to the simpler version found in
\eqref{E:siegel}.

The regularity properties (in the Sobolev sense) of the functions 		$K_\Omega(\cdot,\zeta),\zeta\in\Omega,$ will be determined by the 
pattern of exponents above; details are
given in \S \ref{S:reg} below.

The location and geometry of the complex tangent points are discussed in \S
\ref{S:ctan}.  \S \ref{S:projtran} sets up a change of coordinates serving to reveal
the symmetries of $\Omega$.  In \S \ref{S:Bker} we represent the Bergman
kernel function of (an image of) $\Omega$ as a sum of integrals, and we explain
how the residue calculus can be used to extract the desired asymptotic
expansions, modulo  estimates provided in \S \ref{S:Legendre}.

\medskip

\section{Geometry of complex tangent points}\label{S:ctan}

\medskip
We set
$\rho^2=|a_1|^2+|a_2|^2.$
\medskip

In order for $\partial B_1$ and $\partial B_2$ to intersect non-trivially and real
transversally we must have
\begin{equation}\label{E:rrange}
|1-\rho|<r<1+\rho.
\end{equation}

\medskip

Complex tangent points of $\partial B_1\cap
\partial B_2$ will occur when the radius vectors for  the two spheres are
$\CC$-dependent.   Since the difference of these two vectors is simply the vector
$a$ joining the two centers, it follows that both radii are $\CC$-multiples of
$a$.   In particular,  if $z=(z_1,z_2)$ is a complex tangent point then $z$ is
a multiple of $a$.   Assuming for simplicity that 
$a_1\neq 0$, the complex tangent points of
$\partial B_1\cap
\partial B_2 $ are contained in 
$\left\{\left(z_1, \frac{a_2}{a_1}z_1\right); \;z_1\in \CC\right\}$. 

We shall show that there exist precisely two complex tangent points. Let
$\left(z_1,\frac{a_2}{a_1}z_1\right)\in
\partial B_1\cap \partial B_2$ be a complex tangent point. Then

\begin{align*}
|z_1|^2+\left|\frac{a_2}{a_1}\right|^2 |z_1|^2 &= 1 \\
|z_1-a_1|^2+\left|\frac{a_2}{a_1}z_1-a_2\right|^2 &= r^2 \\
\left|\frac{z_1}{a_1}\right|&= \frac{1}{\rho}\\
\left|\frac{z_1}{a_1}-1\right|&= \frac{r}{\rho}.
\end{align*}

We can write $ \frac{z_1}{a_1}=\frac{1}{\rho}e^{i\chi},\;\; -\pi\le \chi< \pi$. 
From  the last equation of the above system we
obtain

$$
r^2=1+\rho^2-2\rho\, \cos \chi.
$$

\medskip

Since $\cos \chi=\frac{1+\rho^2-r^2}{2\rho}\in(-1,1)$ (by \eqref{E:rrange}), there are
precisely two possible values of $\chi$ in $(-\pi,\pi), \;\chi\neq 0$.  Hence, there exist
precisely two complex tangent points. 

Let  $\theta$ denote the angle between the two radii; thus
$z-a=re^{i\theta}z$.  In particular, $z_1-a_1=re^{i\theta}z_1$ so that
\begin{equation}\label{E:relate}
\frac{e^{i\chi}}{\rho}-1=\frac{z_1}{a_1}-1=re^{i\theta}\frac{z_1}{a_1}
=\frac{r}{\rho}e^{i(\theta+\chi)}
\end{equation}
and 
thus $1-\rho e^{-i\chi}=re^{i\theta}$; taking real and imaginary
parts we find that

\begin{eqnarray}\label{E:chitheta}
1-r \cos \theta&=&\rho\,\cos \chi \notag\\
r \sin \theta &=& \rho\,\sin \chi. 
\end{eqnarray}

Note that the values of $\chi$ and $\theta$ at $q$ can be taken to be the
negatives of the corresponding values at $p$.

The parameters $r$ and $\theta$ have an interpretation extending to more
general situations.  Consider smooth real hypersurfaces $M_1$ and $M_2$ in
$\CC^2$ intersecting real-transversally with a complex tangency at $z\in
M_1\cap M_2$.  Then for suitable $\theta$, the rotation
$R_\theta$ of $T_z \CC^2$ given by multiplication by
$e^{i\theta}$ will map $T_z M_1$ to
$T_z M_2$ and $T_z M_1/H_z$ to $T_z M_2/H_z$, where $H_z=T_z M_j\cap
T_z M_2$ is the maximal complex subspace of both $M_j$.

The Levi-form $\mathcal L_j$ of $M_j$ at $z$ is a hermitian $T_z
M_1/H_z$-valued form on $H_z$; it can be defined by the equation 
\[
\mathcal L_j(X_z)\equiv [X,R_{\pi/2} X]_z \mod H_z
\]
for all smooth vector fields on $M_j$ with values in the maximal complex
subspace of $TM_j$.  Since $\dim_{\mathbb C} H_z=1$, if
$\mathcal L_2$ is non-degenerate then there is $r\in\RR$ so
that
$\mathcal L_1=r R_{-\theta}
\mathcal L_2$.  The parameters $r,\theta$ defined this way match the ones
already defined in the special case of spheres.

% XXXXXXXXXXXXXXXXXXXXXXXXXXXXXXXXXXXXXXXXXXXXXXXXXXXXXXXXXXXXXX
\section{Projective transformation}\label{S:projtran}

\medskip

Let 

\bea 
p:=\left(p_1,\frac{a_2}{a_1}p_1\right),\;\;\;\;\frac{p_1}{a_1}&=&
\frac{1}{\rho}\,e^{i\chi}\\
q:=\left(q_1,\frac{a_2}{a_1}q_1\right),\;\;\;\;\frac{q_1}{a_1}&=&
\frac{1}{\rho}\,e^{-i\chi},
\eea

\medskip
\noindent $-\pi< \chi<\pi,\, \chi\neq 0$ be the two complex tangent points. 
\medskip

We view the balls as embedded in $\CC\PP^2$. If $(z_1:z_2:z_3)$ are the
homogeneous coordinates in $\CC\PP^2$ then the equations of the two balls will
become

\bea
|z_1|^2+|z_2|^2-|z_3|^2 &<& 0\\
|z_1-a_1\,z_3|^2+|z_2-a_2\,z_3|^2-r^2\,|z_3|^2 &<& 0.
\eea

Let $f:\CC\PP^2\to \CC\PP^2$ be a projective transformation such that

\bea
f(0:0:1)&=&\left(q_1:\frac{a_2}{a_1}q_1:1\right)\\
f(0:1:0)&=&\left(p_1:\frac{a_2}{a_1}p_1:1\right).
\eea

\medskip
Let us assume that
$$
f(w_1:w_2:w_3)=(z_1:z_2:z_3)$$

\noindent
where

\begin{eqnarray}\label{E:zwchange}
z_1&=& a_{11}w_1+a_{12}w_2+a_{13}w_3\notag\\
z_2&=&a_{21}w_1+a_{22}w_2+a_{23}w_3\\
z_3&=& a_{31}w_1+a_{32}w_2+a_{33}w_3.\notag
\end{eqnarray}

\smallskip 
Due to the above constraints we see that the matrix $A$ of the transformation $f$ is

\[ \left(
\begin{array}{ccc}
a_{11} & \lambda p_1 & \mu q_1 \\ a_{21} & \lambda \frac{a_2}{a_1}p_1 & \mu \frac{a_2}{a_1}q_1\\ a_{31} & \lambda & \mu 
\end{array}\right).  \]

\medskip
The equation of the first ball shall be transformed under $f$ to

\begin{align*}
(|a_{11}|^2+|a_{21}|^2&-|a_{31}|^2)|w_1|^2\\
&+ 2\,{\Re}\left[\overline\lambda\left\{\overline p_1\left(a_{11}+a_{21}
\overline{\left(\frac{a_2}{a_1}\right)}\right)-a_{31}\right\} w_1\,\overline
w_2\right]\\ 
&+ 2\,{\Re}\left[\mu\left\{q_1\left(\overline a_{11}+\overline
a_{21}\frac{a_2}{a_1}\right) -\overline a_{31}\right\}w_3\,\overline
w_1\right]\\
&+2\,{\Re}\left[\lambda\,\overline\mu \left(p_1\overline
q_1+\left|\frac{a_2}{a_1}\right|^2 p_1\overline q_1-1\right)w_2\overline
w_3\right]\\ &<0.
\end{align*}

\medskip
If we require that the $w_1$ direction  be 
tangent to $\partial B_1$ at $(0:1:0)$
and at $(0:0:1)$ then 

\bea
\overline\lambda\left\{\overline p_1\left(a_{11}+a_{21}
\overline{\left(\frac{a_2}{a_1}\right)}\right)-a_{31}\right\}&=& 0\\
\mu\left\{q_1\left(\overline a_{11}+\overline a_{21}\frac{a_2}{a_1}\right)
-\overline a_{31}\right\}&=&0.
\eea

Since $\lambda,\mu \in \CC^*, \;p_1\neq q_1$ the above system yields

\bea
a_{31}&=&0\\
a_{11}&=& -a_{21}\overline{\left(\frac{a_2}{a_1}\right)}.
\eea

\medskip

The matrix $A$ of the transformation $f$ shall become

\begin{equation}\label{E:matr}
\left(
\begin{array}{ccc}
-\overline{\left(\frac{a_2}{a_1}\right)}a_{21} & \lambda p_1 & \mu q_1 \\
a_{21} &
\lambda \frac{a_2}{a_1}p_1 & \mu \frac{a_2}{a_1}q_1\\ 0 & \lambda & \mu 
\end{array}\right) 
\end{equation}

\medskip 
\noindent
with $\lambda,\mu \in \CC^*$. 

\medskip 

We shall normalize the coefficient $|a_{11}|^2+|a_{21}|^2-|a_{31}|^2=
\frac{\rho^2}{|a_1|^2}\,|a_{21}|^2$ of $|w_1|^2$  such that it equals $1$.
This will imply that

\begin{equation}\label{E:a21cond}
|a_{21}|=\frac{|a_1|}{\rho}.
\end{equation}

\medskip

We shall also choose $\lambda,\mu\in \CC^*$ such that
$$
\lambda\,\overline{\mu}\left\{p_1\,\overline q_1
\left(1+\left|\frac{a_2}{a_1}\right|^2\right)-1\right\}=-1
$$

\medskip 
\noindent
or equivalently (using the fact that $\frac{p_1}{a_1}=\frac{1}{\rho}\,e^{i\chi},\; \frac{q_1}{a_1}=\frac{1}{\rho}\,e^{-i\chi}$)

\begin{equation}\label{E:lammucond}
\lambda\overline{\mu}( e^{2i\chi}-1)=-1.
\end{equation}

\medskip 
\noindent
We can choose 
\begin{equation}\label{E:lammudef}
\lambda=\frac{1}{ e^{i\chi}-1}=\frac{ e^{-i \frac{\chi}{2}}}{ 2 i \sin
\frac{\chi}{2}},\qquad
\mu=\frac{-1}{e^{-i \chi}+1}=\frac{-e^{i
\frac{\chi}{2}}}{2 \cos \frac{\chi}{2}}.
\end{equation}

{\it Note:} This is possible since $\left|\cos\chi\right|<1$. 

\medskip

Thus, the first ball is now described by the equation
\begin{equation}\label{E:sieg1}
|w_1|^2-2\,\Re w_2\overline{w_3}<0.
\end{equation}

\medskip
\noindent The equation of the second ball is transformed under $f$ to 
\begin{align*}
|a_{21}|^2 &\left(1+\left|\frac{a_2}{a_1}\right|^2\right)|w_1|^2\\
&+
2\,\Re\left[\lambda\overline{\mu}\left\{(p_1-a_1)\overline{(q_1-a_1)}\left(1+
\left|\frac{a_2}{a_1}\right|^2\right)-r^2\right\}w_2\overline{w_3}\right]\\ 
&< 0.
\end{align*}

\medskip
From the previous
normalization we have that $|a_{21}|^2 
(1+ |\frac{a_2}{a_1} |^2 )=1$, so
the equation of the second ball becomes
\begin{equation}\label{E:transball}
|w_1|^2+2\Re\phi w_2\overline{w_3}<0
\end{equation}
with
\bea
\phi&=&\lambda\overline{\mu}\left\{(p_1-a_1)\overline{(q_1-a_1)}\left(1+
\left|\frac{a_2}{a_1}\right|^2\right)-r^2\right\}\\
&=&
\lambda\overline{\mu}\left\{\rho^2
\left(\frac{p_1}{a_1}-1\right)\overline{\left(\frac{q_1}{a_1}-1\right)}
-r^2\right\}.
\eea

\medskip

Recalling from \eqref{E:relate} that 

\bea
\frac{p_1}{a_1}-1&=& \frac{r}{\rho}e^{i(\theta+\chi)}\\
\frac{q_1}{a_1}-1&=& \frac{r}{\rho}e^{-i(\theta+\chi)}
\eea

\medskip
\noindent and using 
\eqref{E:lammucond},
we can rewrite $\phi$  as

$$
\phi=-r^2\frac{e^{2i(\theta+\chi)}-1}{e^{2i\chi}-1}
=-r^2e^{i\theta}\frac{\sin(\theta+\chi)}{\sin\chi}.
$$
\medskip

Using the angle addition formula and the identities in
\eqref{E:chitheta} we find that

\begin{equation*}\label{E:phiform}
\phi=-re^{i\theta}.
\end{equation*}

Let $B_1', B_2'$ be the preimages of $B_1, B_2$ under 
the projective transformation $f$ and standard normalizations. Let
$\Omega'=B_1'\cap B_2'$. Since $w_3\neq 0$ in $\Omega'$  (in view of
\eqref{E:sieg1} and 
\eqref{E:transball}), setting 
$\tilde w_1=\frac{w_1}{w_3},\;\; \tilde w_2=\frac{w_2}{w_3}$   we  may
summarize our work as follows.

\medskip

\begin{proposition}
The inverse of the  projective transformation induced 
by the matrix \eqref{E:matr} subject to \eqref{E:a21cond} and 
\eqref{E:lammudef} maps $\Omega$
to 
$$
\hat\Omega'={\mbox{affine}}\; 
\Omega'=\{(\tilde w_1, \tilde w_2); \;|\tilde w_1|^2<{\min}\{2\Re\tilde
w_2,\;-2\Re\,
\phi  \tilde w_2 \}\}.
$$
\end{proposition}

\begin{proof}[Proof of Theorem \ref{T:maingen}, part (i).]
The domain $\hat\Omega'$ admits the $\RR\times S^1$ action
\begin{equation*}
(s,\theta)\cdot(\tilde w_1,\tilde w_2)=(s e^{i\theta}\tilde w_1,s^2\tilde w_2).
\end{equation*}
This action pulls back to a $\RR\times S^1$ action on $\Omega$.

The result now follows by application of [Ba1, proof of Theorem 3]
to this action.  (It also follows from later computations in this paper.)
\end{proof}

For future reference we note that the inverse of the map \eqref{E:zwchange} is
given by 
\begin{align}\label{E:inverse}
w_1&=-\frac{z_1 a_2-z_2 a_1}{ a_1 a_{21} (1+|\frac{a_2}{a_1}|^2)}\notag\\
w_2&= \frac{ z_1+\overline{(\frac{a_2}{a_1})}z_2-q_1 
(1+|\frac{a_2}{a_1}|^2) z_3}{ \lambda (p_1-q_1) (1+|\frac{a_2}{a_1}|^2)}\\
w_3&= \frac{ z_1+\overline{(\frac{a_2}{a_1})} z_2-p_1
( 1+|\frac{a_2}{a_1}|^2) z_3}{ -\mu (p_1-q_1)(1+|\frac{a_2}{a_1}|^2)}\notag.
\end{align}

\section{The Bergman kernel of the intersection of two balls}\label{S:Bker}

Let $\Phi:\hat \Omega'\to\CC^2$ be the transformation defined by $\Phi(\tilde 
w_1, \tilde w_2)=(t_1,t_2)$ with

\bea
t_1&=& \frac{1}{\sqrt{2}}\,\tilde w_1\,\tilde w_2^{-\frac{1}{2}},\;\; t_2=u+iv=\log\tilde w_2 \\
\tilde w_1&=& \sqrt{2}e^{\frac{t_2}{2}}\,t_1,\;\; \tilde w_2=e^{t_2}.
\eea

\medskip
Then $D':=\Phi(\hat \Omega')$ is defined by the inequality
\begin{equation*}
\left\{|t_1\right|^2<
\psi_{r,\theta}(v)\}.
\end{equation*}

\medskip
\noindent where
\;$\psi_{r,\theta} (v)={\min}\{\cos v, r\cos (v+\theta)\}$.

$D'$ is a Hartogs domain invariant under the rotations $(t_1, t_2)\mapsto ( e^{i
\alpha} t_1, t_2)$. By Fourier expansion the Bergman space $\mathcal{ H}(D')$ 
(the space consisting of all square-integrable, holomorphic functions
in $D'$)
admits  an orthogonal decomposition 

$$
\mathcal{H}(D')=\oplus \,\mathcal{H}_j( D')
$$

\noindent where $\mathcal{H}_{j}(D')$ is the subspace consisting of all
square-integrable, holomorphic functions $f$ in $D'$ that satisfy $f(
e^{i\alpha} t_1, t_2)=e^{i j \alpha} f(t_1,t_2)$. Functions with this 
property are of the form $f(t_1, t_2)=t_1^{j} f_{1}(t_2)$, \; $f_{1}$ 
holomorphic in $t_2$.  The Bergman kernel
$K_{D'}(t,\tau)$ satisfies
$$
K_{D'}(t,\tau)=\sum_{j\ge 0} K_{j}(t,\tau),
$$

\noindent where $K_j (t,\tau)$ is the reproducing kernel for $\mathcal{H}_j
(D')$. 

\medskip

\medskip

Using an argument similar to the one in Section 1 of [Ba2] and noting  that for
$f,g$ holomorphic functions in $t_2$ we have

\bigskip

\begin{multline*}
\int_{D'} f(t_2)t_1^j\, \overline{g(t_2)t_1^{k}}\,dV\\
=
\begin{cases}
0 & j\neq k,\\
\frac{\pi}{j+1} \int\limits_{v_{{\min}}<v<v_{{\max}}} {f(t_2)\overline{g(t_2)}
\psi_{r,\theta}^{j+1}(v)\,dA} & j=k,
\end{cases}
\end{multline*}
 we find
that

\begin{align}\label{E:D'formula}
K_{D'}&((t_1,t_2),(\tau_1,\tau_2))\notag\\
&=
\frac{1}{2\pi^2}\sum_{j\ge 0}t_1^j\,\overline{\tau_1}^j \,(j+1)
\int\limits^{\infty}_{-\infty} \frac{e^{i(t_2-\overline{\tau}_2)\, \xi}\;d\xi}
{\int\limits_{v_{{\min}}<v<v_{{\max}}} \psi_{r,\theta}^{j+1}(v)\,e^{-2v\xi}
\,dv}\notag\\
&=\frac{1}{4\pi^2}\sum_{j\ge 0}t_1^j\,\overline{\tau_1}^j \,(j+1)
\int\limits^{\infty}_{-\infty} \frac{e^{i(t_2-\overline{\tau}_2)\,
\frac{\xi}{2}}\;d\xi} {\int\limits_{v_{{\min}}<v<v_{{\max}}}
\psi_{r,\theta}^{j+1}(v)\,e^{-v\xi}
\,dv}.
\end{align}

\medskip

To simplify
notation we shall write  from now on $( v_{{\min}}, v_{{\max}}):=J$ and
$\eta_{j+1}(v):=-(j+1)\log\psi_{r,\theta}(v)$.

\medskip

Let us assume for the moment that we can apply contour integration arguments 
to each one of the above integrals for appropriate $t_2,\tau_2$.  Then for $h>0$
we have:

\begin{eqnarray}\label{E:D'asymp}
\int\limits^{\infty}_{-\infty} \frac{e^{i(t_2-\overline{\tau}_2)\,
\frac{\xi}{2}}\;d\xi} {\int\limits_{J} e^{ (-v\xi-\eta_{j+1}(v))} \,dv}&=&
-2\pi i  \sum_{-h< \Im\xi<0 }
\Res\left(\;\frac{e^{i(t_2-\overline{\tau}_2)(\frac{.}{2})}}{F_{j+1}(.)},\; \xi\;\;
\right) 
 +\notag\\
&+&\int\limits^{\infty}_{-\infty} \frac{e^{i(t_2-\overline{\tau}_2)\, 
\frac{(-x-ih)}{2}}\;dx}{F_{j+1}(-x-i h)} 
\end{eqnarray}

\medskip

\noindent where 

\begin{equation}\label{E:Fjdef}
F_{j+1}(\xi):=\int\limits_{J} e^{-v\xi-\eta_{j+1}(v)} dv
\end{equation}
$\left(\text{hence }
F_{j+1}(-x-i h):=\int\limits_{J} e^{iv h} e^{v x-\eta_{j+1}(v)} dv\right).
$

\medskip

We will see below in Corollary \ref{C:finstrip} that the union of the zero sets 
of the
$F_{j+1}$ is finite (counting multiplicity) in each strip
$-h\le \Im\xi\le 0$.

In particular, for all but a discrete set of $h$, we have
\begin{equation}\label{E:wcond}
F_{j+1}(-x-ih)\ne0 \text{ for all } x\in {\mathbb R},\;  j\ge 0
\end{equation}
so that the final
integrand in
\eqref{E:D'asymp} does not suffer a vanishing denominator.

\medskip

In \S \ref{S:Legendre} we shall show that the use of the residue theorem above is
valid,  and in Proposition \ref{P:sum} we show that when the residue expansions
\eqref{E:D'asymp} are substituted into \eqref{E:D'formula}, the  sum of 
integrals is an  error term of magnitude $ O (
e^{ t_2 h/2})$ as $\Re t_2\to-\infty$, uniformly as $\tau$ ranges over any compact
subset of $D'$.

\medskip 

Using the last remark and applying the transformation formula 

\begin{equation}\label{E:Btran}
K_{\hat\Omega'}(\tilde w, \tilde{\omega})=\frac{1}{2 ( \tilde w_2 \overline{
\tilde\omega_2})^{\frac{3}{2}}} K_{D'} ( \Phi(\tilde w), \Phi(\tilde\omega))
\end{equation}

\noindent for the Bergman kernel we obtain:

\begin{multline*}
K_{\hat \Omega'}((\tilde w_1,\tilde w_2),(\tilde \omega_1, \tilde \omega_2))
\\
=\tfrac{1}{4\pi i\;(\tilde w_2\,\overline{\tilde
\omega_2})^{\frac{3}{2}}}\;
\sum_{j\ge 0}
\sum_{-h<\Im \xi<0}
\tfrac{j+1}{2^j} 
\tfrac{(\tilde w_1\,\overline{\tilde\omega_1})^{j}} {(\tilde
w_2\,\overline{\tilde
\omega_2})^{\frac{j}{2}}}\; \Res\left( \tfrac{e^{i( \log \tilde w_2-\log
\overline {\tilde\omega}_2)(\frac{.}{2})}} {F_{j+1}(.)}, \xi\right)+\\
+O\left(
\tilde w_2^{\frac{h-3}{2}}
\right)
\end{multline*}
as $\tilde w\to 0$ in
$\hat\Omega'$, uniformly as $\tilde\omega$ ranges over any compact subset of $\hat
\Omega'$.

If the zeroes of the $F_{j+1}$ in the strip $-h<\Im \xi<0$ are all simple, the
expansion may be written in the form
\begin{multline*}
K_{\hat \Omega'}((\tilde w_1,\tilde w_2),(\tilde \omega_1, \tilde \omega_2))
\\
=\frac{1}{(\tilde w_2\,
\overline{\tilde \omega_2})^{\frac{3}{2}}}\;
\sum_{j\ge 0}
\sum_{k}
c_{j,k}  \frac{(\tilde
w_1\,\overline{\tilde\omega_1})^{j}} {(\tilde w_2\,\overline{\tilde
\omega_2})^{\frac{j}{2}}}\;
\; \left(\frac{\tilde w_2}{\overline{\tilde\omega}_2}\right)^{i
\frac{\xi_{j,k}}{2}}+\;
\\
+O\left(
\tilde w_2^{\frac{h-3}{2}}
\right),
\end{multline*}

\medskip
\noindent where $\xi_{j,k}$ is an enumeration of zeros in the strip
$-h<\Im\xi<0$, and the
$c_{j,k}$ are constants that arise from the residue principle. 

In the general case, if $F_{j+1}(\xi)$ has a root of multiplicity $m_k$ at
$\xi_{j,k}$ then our expansion looks like:

\begin{multline}\label{E:siegel}
K_{\hat \Omega'}((\tilde w_1,\tilde w_2),(\tilde \omega_1, \tilde \omega_2))
\\
=\frac{1}{(\tilde w_2\,
\overline{\tilde \omega_2})^{\frac{3}{2}}}\;
\sum_{j\ge 0}
\sum_{k}
\frac{(\tilde
w_1\,\overline{\tilde\omega_1})^{j}} {(\tilde w_2\,\overline{\tilde
\omega_2})^{\frac{j}{2}}}\;
\; \left(\frac{\tilde w_2}{\overline{\tilde\omega}_2}\right)^{i
\frac{\xi_{j,k}}{2}}\;P_{j,k}\left(\log\frac{\tilde
w_2}{\overline{\tilde\omega_2}}\right)+
\\
+O\left(
\tilde w_2^{\frac{h-3}{2}}
\right),
\end{multline}
where $P_{j,k}$ is a polynomial of degree $m_k-1$.

\medskip

We are only a step away from our original goal of obtaining the asymptotic
expansion for  the Bergman kernel on the intersection of the two balls in
$\CC^2$.  Using \eqref{E:inverse} and recalling  that 

$$
\tilde w_1= \frac{w_1}{w_3},\;\;\;\; \tilde w_2=\frac{w_2}{w_3}
$$

\noindent we can rewrite $\tilde w_1,\; \tilde w_2$ as functions of
$z_1,z_2,z_3$.  Setting $z_3=1$ we can obtain a  biholomorphic map  

\bea\label{E:Psidef}
\Psi: B_1\cap B_2 &\to& \hat\Omega' \\
(z_1,z_2) &\mapsto& (\tilde w_1, \tilde w_2)\\
\eea

\medskip
\noindent
where
\bea
\tilde w_1&=&\frac{2 i r \mu \; \sin \theta  \;\;< z-q, ^{T}a>}{  a_1 a_{21}(1+|\frac{a_2}{a_1}|^2)\;\;(<z-q,a>-2i r \sin \theta) }\\
\tilde w_2&=& -\frac{ \mu \;\; <z-q,a> }{ \lambda \;\;( <z-q, a>-2i r \sin
\theta)}.\\
\eea

\medskip
\noindent     
Applying the transformation formula for the Bergman kernel we obtain our main
result.

\begin{theo}\label{T:main}  Let $h>0$ satisfy \eqref{E:wcond}.  Then

\begin{multline}\label{E:fullstory}
K_{\Omega}(z,\zeta)\\
=
\sum_{j}\sum_{-h<\Im \xi<0}
\frac{< z-q, ^{T} a> ^{j}\;\;  \overline{<\zeta-q,
^{T}a>}^{j}}{\;\;< z-q,  a> ^{-\frac{i \xi_{j,k}}{2}+\frac{j+3}{2}}\;\;  \overline{<\zeta-q,
a>}^{\frac{i \xi_{j,k}}{2}+\frac{j+3}{2}}}\cdot\\
\cdot
\left(< z-q,  a>- 2i r \sin \theta \right)
^{-\frac{i\xi_{j,k}}{2}-\frac{j+3}{2}}\cdot\\
\cdot
\overline{\left(< \zeta-q,  a> -2i r \sin \theta \right)}
^{\frac{i\xi_{j,k}}{2}-\frac{j+3}{2}}\cdot\\
\cdot
P_{j,k}\left(\log
\frac{< z-q, a>\; \overline{\left(<\zeta-q,a>-2 i r \sin \theta \right)}}
{\overline{< \zeta-q,  a>}\;\; \left(<z-q,  a>-2 i r \sin \theta \right)}
\right)\\
+O\left(<z-q,a>^{\frac{h-3}{2}}\right)
\end{multline}
as $z\to q$, uniformly as $\zeta$ ranges over any compact subset of $\Omega$;
here  
\begin{itemize}
\item $q$ is a complex tangent point of the intersection of the two balls;
\item $a=(a_1,a_2)$;
\item $^{T}a:=(\overline{a_2}, -\overline{a_1})$;
\item $\{\xi_{j,k}\}$  is an enumeration of the zeroes of $F_{j+1}$
in the strip $-h<\Im\xi<0$;
\end{itemize}
and
\begin{itemize}
\item $P_{j,k}$ is a polynomial of degree one less than the multiplicity of
$F_{j+1}$ at $\xi_{j,k}$.
\end{itemize}
\end{theo}

\medskip
For brevity we set
\begin{equation}\label{E:Jdef}
J:=(\upsilon_{\min}, \upsilon_{\max})=
\begin{cases}
(-\frac{\pi}{2},
\frac{\pi}{2}-\theta), &\theta>0\\
(-\frac{\pi}{2}-\theta,
\frac{\pi}{2}), &\theta<0.
\end{cases}
\end{equation}
\medskip
\begin{lem}\label{L:goodstrip}
$F_{j+1}(\xi)$ does not vanish in the strip $\left|\Im\xi\right| \le\frac{\pi}{|J|}$.
\end{lem}

\begin{proof}
Let $v^*$ denote the midpoint of $J$.  Then 
\begin{align*}
\Re e^{v^*\xi} F_{j+1}(\xi)&= \int_J \Re e^{(v^*-v)\xi-\eta_{j+1}(v)}\,dv
\\&>0
\end{align*}
(since the integrand is positive).
\end{proof}

\medskip

\begin{proof}[Proof of Theorem \ref{T:maingen}.]
The expansion \eqref{E:crudeasym} is obtained by expanding the non-logarithmic
factors in  \eqref{E:fullstory} in powers of $\langle z-q,
^{T}a\rangle$ or $\langle z-q, a\rangle$.  The assertion about the location of
the $\gamma_j$ follows from Lemma \ref{L:goodstrip}.
\end{proof}

If we make a change of coordinates that  sends 

\bea
(0,0) &\mapsto& (q_1,q_2) \\
(1,0) &\mapsto&  (p_1,p_2)\\
\eea

\noindent then in the new coordinates the Bergman kernel will look like

\begin{multline*}
K((\zeta_1, \zeta_2),(\zeta'_1, \zeta'_2))\sim\\ \tfrac{1}{ (\zeta_1 (\zeta_1-1)
\overline{\zeta'_1}(\overline{\zeta'_1}-1))^{\frac{3}{2}}} \sum_{j,k} 
\left(
\tfrac{{\zeta_2}^2 \overline{\zeta'_2}^2}{(\zeta_1 (\zeta_1-1)
\overline{\zeta'_1}(\overline{\zeta'_1}-1))}\right)^{\frac{j}{2}}
\left(\tfrac{\zeta_1}{\zeta_1-1}
\tfrac{\overline{\zeta'_1}-1}{\overline{\zeta'_1}}\right)^{\tfrac{ i\xi_{j,k}}{2}} \;
P_{j,k} \left(\log \tfrac{ \zeta_1 \;(\overline{\zeta'_1}-1)}{ (\zeta_1-1)\;
\overline{\zeta'_{1}}}\right) 
\end{multline*}

\noindent in a neighborhood of $(0,0)$.

The results from [Ba3] can be used to provide a more operator-theoretic approach
to obtaining such asymptotic expansions.

\medskip

We conclude this section with some remarks on the location of the zeroes
$\xi_{j,k}$ in special cases.

\medskip

First, in the special case (not otherwise allowed
in this paper) $\theta=0$ where $\Omega$ is a ball we have
\begin{equation*}
F_{j+1}(\xi)=
\begin{cases}
\dfrac{2(j+1)!\cosh(\pi\xi/2)}{ (1+\xi^2) (9+\xi^2)\cdots
((j+1)^2+\xi^2)},
& j+1\text{ odd};\\
\,\\
\dfrac{2(j+1)!\sinh(\pi\xi/2)}{
 \xi (4+\xi^2)(16+\xi^2)\cdots ((j+1)^2+\xi^2)},
& j+1\text{ even}.
\end{cases}
\end{equation*}
\medskip

In particular, the zeroes of each $F_{j+1}$ are simple.

It follows by use of the argument principle that in any strip $-h\le\Im\xi\le 0$, all
the 
$F_{j+1}$ have simple zeroes provided that $|\theta|$ does not
exceed $\theta_0(h)$.

\medskip

Let us now consider the case $r=1$, $\theta\ne0$.   In this case we have
\begin{equation*}
F_1(\xi)=
\begin{cases}
\frac{
\left(1+ie^{-\frac{\pi-\theta}{2}\xi+i\frac{\theta}{2}}\right)
\left(1-ie^{-\frac{\pi-\theta}{2}\xi-i\frac{\theta}{2}}\right)
}
{1+\xi^2}
\cdot e^\frac{\pi\xi}{2}, & \theta>0\\
\frac{
\left(1+ie^{\frac{\pi+\theta}{2}\xi-i\frac{\theta}{2}}\right)
\left(1-ie^{\frac{\pi+\theta}{2}\xi+i\frac{\theta}{2}}\right)
}
{1+\xi^2}
\cdot e^\frac{-\pi\xi}{2}, & \theta<0.
\end{cases}
\end{equation*}

\noindent 
The zeroes of $F_1$ are given by
\begin{equation}\label{E:zs}
i \xi_{0,n}=
\begin{cases}
\left( - 1+\frac{2 (n+1) \pi}{\pi-|\theta|}\right), &n \text{ odd;}\\
\left( 1+\frac{2 n \pi}{\pi-|\theta|}\right), &n \text{ even.}
\end{cases}
\end{equation}

\noindent Using repeated integration by parts, the general form of the $F_{j+1}$ for
 $r=1$ can
be given as follows.

\begin{enumerate}
\item[(i)]  If $j+1$ is odd, then 

$$
F_{j+1}(\xi)= C e^{\theta \frac{\xi}{2}} \frac{\cosh( (\pi-\theta)
\frac{\xi}{2})-P_{j}^\theta(\xi)} { (1+\xi^2) (9+\xi^2)\cdots
((j+1)^2+\xi^2)}
$$

\medskip
\noindent where $P^\theta_{j}(\xi)$ is an even polynomial  in $\xi,$ of degree $j,$ 
whose coefficients depend on $\theta$ and $C$ is  some absolute constant.  

\medskip

\item[(ii)]  If $j+1$ is even, then 

$$
F_{j+1}(\xi)=C' e^{(\pi-\theta)\frac{\xi}{2}} 
\frac{\sinh( (\pi-\theta) \frac{\xi}{2})-Q^\theta_{j}(\xi)}{
 \xi (4+\xi^2)(16+\xi^2)\cdots ((j+1)^2+\xi^2)}
$$

\medskip
\noindent  where $Q_j^\theta (\xi)$ is an odd polynomial in $\xi,$ of degree
$j,$ whose  coefficients depend on $\theta$ and $C'$ is  some absolute
constant.    

\end{enumerate}

\medskip

The above formulas indicate that the zeroes lying off the imaginary
axis will occur in pairs $\xi, -\overline\xi$; the pair will always lead to
terms of equal strength (as far as it concerns estimates) so both terms should  be
considered together.

We cannot provide an explicit formula for the first zero of
$F_2$, but for $r=1, \theta\sim 0$   the first zero can be approximated by
\begin{equation*}
i\xi_{1,1}=4+\frac{8}{\pi}
|\theta|+O(\theta^{2}).
\end{equation*}
Similarly, for $r=1, \theta=\pi-\epsilon,
\epsilon\sim 0^+$, we find by rescaling $\xi$  that the conjugate pair $i\xi_{1,1},
i\xi_{1,2}$ can be approximated by 
\begin{equation*}
\dfrac{ 14.995 \ldots\pm i
5.537\ldots}{\epsilon}+O(1).
\end{equation*}

\bigskip

%Mac graphics commands
%\input picmacs.tex
%
%\centerpicture 6.26in by 3.86in (roots scaled 600)

\begin{figure}
\scalebox{.70}
{\includegraphics{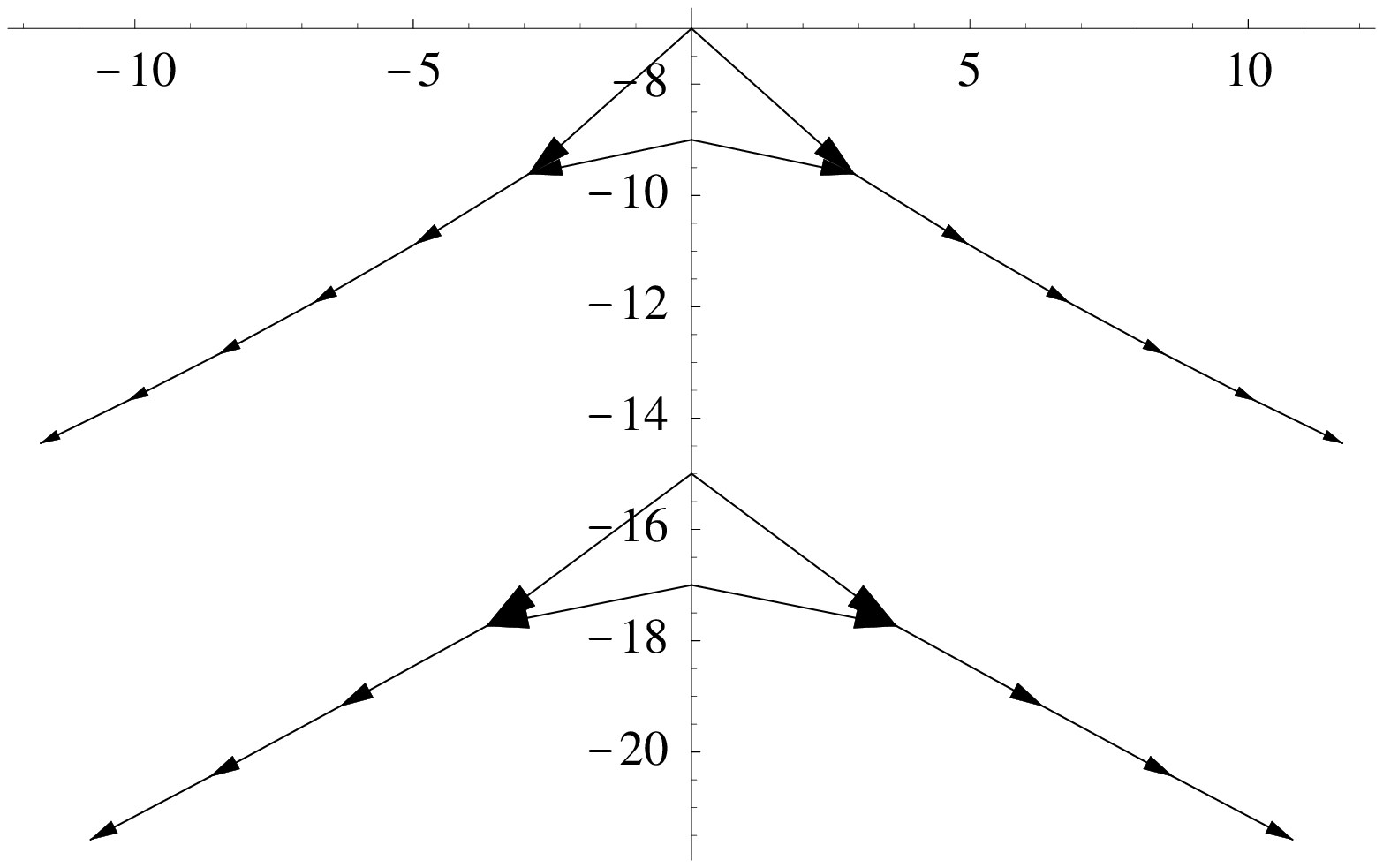}} 
\caption{}
\end{figure}

For general values of $\theta$ the roots can at least be explored numerically.

Figure 1 shows a portion of the root pattern for $r=1, \theta=\frac{\pi}{2}$.
The arrows connect zeroes of $F_j$ to  nearby zeroes of $F_{j+1}$.  The  zeroes of $F_1$ lie on the
imaginary axis at $(8m\pm 1)i$.  Passing to $F_2$, the first four zeroes (at least) are converted into two
root pairs, symmetric across the imaginary axis.  As $j$ increases further, the roots move down and out.

\medskip

When $r>1$ we can
use as before repeated integration by  parts to compute $F_{j+1}$ -- however
there won't be such nice formulas  as above. $F_{j+1}$ can be written as a finite
sum of exponential functions
$$
\sum_{j} e^{ a_{j} \frac{\xi}{2}} P_{j}(\xi)
$$
\medskip
where the $P_{j}(\xi)$ are rational functions in $\xi$, whose coefficients
depend  on $\theta, r$.

\medskip
When $r>1,\;\theta\sim 0$ we can get the following information

\bea
i \xi_{0,k}&=& 1+ 2 k +i\frac{ 2  r
k(k+1)}{(r-1)\pi}
\,\theta|\theta|+O(\theta^{3})\\
i \xi_{1,k} &=& 2 +2k +O(\theta^{3}),\\
\eea

\medskip
\noindent while when $r>1,\; \theta=\pi-\epsilon, \; \epsilon\sim 0^+$ we can
rescale as before to obtain 

\bea
i \xi_{0,k}&=& \frac{a_{k}}{\epsilon}+O(1)\\
i \xi_{1,k} &=& \frac{b_{k}}{\epsilon}+O(1).\\
\eea

\medskip

For $\theta$ close enough to zero, Proposition \ref{P:est} below can be combined
with a Rouch\'e's theorem argument to show that $\xi_{0,1}$ will have the
largest imaginary part among roots in the lower half-plane and thus will
provide the key to regularity properties of the kernel function. (See \S
\ref{S:reg} below.)

% XXXXXXXXXXXXXXXXXXXXXXXXXXXXXXXXXXXXXXXXXXXXXXX

\medskip
\section{Estimates for Fourier-Laplace transforms of log-convex functions}
\label{S:Legendre}

We would like to use Residue Calculus to obtain asymptotic expansions
for every 
$j$ for the
integrals

\begin{equation}\label{E:keyint}
\int\limits^{\infty}_{-\infty} \frac{e^{i(t_2-\overline{\tau}_2)\,
\frac{\xi}{2}}\;d\xi} {F_{j+1}(\xi)};
\end{equation}

\medskip

\noindent here $F_{j+1}$ is the function defined in \eqref{E:Fjdef}, in which 
$\eta_{j+1}(v)$ is the  piecewise
$C^2$, strictly  convex function $-(j+1)\;\log\psi_{r,\theta}(v)$.

For convenience we will focus on the case $\theta>0$.  Also, we assume that
$(t_1,t_2)\in \overline{D'}$, $(\tau_1,\tau_2)\in {D'}$; it follows that
 \begin{equation}\label{E:Imtbnd}
\Im(t_2-\overline\tau_2)\in(-\pi,\pi-2\theta).
\end{equation} 

We are interested in the
behavior of the above integrals  as $\Re t_2\to-\infty$.  
To get started we shall need lower bounds for 
$
|F_{j+1}(-x-i h)|:=\left|{\int\limits_{J} 
e^{iv h} e^{vx-\eta_{j+1}(v)}\,dv} \right|.
$

\medskip

Recall that

\[ \eta_{j+1}(v)= \left\{ \begin{array}{ll}
-(j+1)\log\cos v   & \;\;\;\;  -\frac{\pi}{2} <  v \le v_{0}\\
-(j+1)\log[r\cos(v+\theta)] & \;\; \;\; v_{0} \le v < \frac{\pi}{2}-
\theta  \\
\end{array}
\right. \]

\medskip

\noindent where $v_{0}=\arctan\left[\frac{r\cos \theta-1}{ r \sin
\theta}\right],$ i.e. $v_{0}$ is the point in $J$ where $\cos v $ and
 $r\cos(v+\theta)\;$ intersect. 

For convenience, we set $\eta_{j+1}(v)=+\infty$ for $v\notin J$.

\medskip

We are going to use the Legendre transform 
\begin{equation*}\label{E:Legdef}
\widetilde {\eta_{j+1}}(\xi):={\max} \{ v\xi-\eta_{j+1}(v); \;v\in J \}
\end{equation*}
of $\eta_{j+1}$. 

\medskip

At the points where $\eta_{j+1}$ is differentiable we can compute the Legendre 
transform using differential calculus (see [H\"o], pages 16-19). 

\medskip

We have:
\begin{multline}\label{E:Legt}
\widetilde{\eta_{j+1}}(\xi)=\\
=\begin{cases}
\xi \arctan \frac{\xi}{j+1}-\frac{j+1}{2} \log (1+(\frac{\xi}{j+1})^{2})\\
\hspace{1.0in}\text{ for } \xi< (j+1)\tan v_{0};\\
v_{0} \xi+(j+1) \log \cos v_{0}\\
\hspace{1.0in}\text{ for }(j+1)\tan v_{0} \le  \xi  \le (j+1) \tan (v_{0}+\theta);\\
\xi (\arctan \frac{\xi}{j+1}-\theta)-\frac{j+1}{2} \log (1+
(\frac{\xi}{j+1})^{2})+\log  r^{j+1}\\
\hspace{1.0in}\text{ for }
 \xi>(j+1)\tan (v_{0}+\theta).
\end{cases}
\end{multline}

\medskip
\noindent

Since $\eta_{j+1}$ is strictly convex, there is a unique $\mu_{j+1}(\xi)$ satisfying

$$
\widetilde{\eta_{j+1}}(\xi)=\xi \mu_{j+1}(\xi)-\eta_{j+1}(\mu_{j+1}(\xi));
$$

\medskip 
\noindent in fact, 
\begin{equation}\label{E:muspec}
\mu_{j+1}(\xi)=\left\{ \begin{array}{ll}
\arctan \frac{\xi}{j+1}&  \phantom{sadsdasdasdasa}\xi< (j+1) \tan v_{0}\\
v_{0}& (j+1)\tan v_{0}\le \xi \le (j+1) \tan (v_{0}+\theta)\\
\arctan \frac{\xi}{j+1}-\theta & \phantom{sadasdasdasdad}\xi>(j+1) \tan (v_{0}+\theta).
\end{array}
\right.
\end{equation}

\medskip

We shall show the following:

\begin{proposition} \label{P:est} Let 
$\widetilde{\eta_{j+1}}:={\max} \{ v\xi-\eta(v); \;v\in J \}$ denote the
Legendre transform of $\eta_{j+1}$. Then
\begin{enumerate}
\item[i)] For every $h_0>0$,  there exists a positive constant $L=L(h_0)$ such that for 
$j\ge
j_0(h_0,|J|)$ sufficiently large,
$0\le h \le h_0$, and 
$x<(j+1)\;\tan v_{0} $ or $x>(j+1)\;\tan (v_{0}+\theta)$, we have 

\begin{equation}\label{E:noncorner}
|F_{j+1}(-x-ih)|\ge \frac{L}{\sqrt{  \eta_{j+1}''(\mu_{j+1}(x))}}\;
e^{\,\tilde\eta_{j+1}(x)}.%%%%e^{\,\widetilde{\eta_{j+1}}(x)}.
\end{equation}

\item[ii)]  For every $h_0>0$, there exists a positive constant $K=K(h_0)$ such that for $j\ge
j_0(h_0,|J|)$ sufficiently large, $0\le h \le h_0$, and $(j+1)\;\tan v_{0}\le x\le
(j+1)\;\tan (v_{0}+\theta)$, we have:

$$
|F_{j+1}(-x-i h)|\ge \frac{K}{ j+1}\; 
e^{\,\tilde\eta_{j+1}(x)}.%%%%e^{\,\widetilde{\eta_{j+1}}(x)}.
$$
\end{enumerate}

\end{proposition}

\medskip

{\bf Remark:} {\em For any $j$, $i)$ still holds if we take $ |x|> x_{0}(j,h_0,|J|)$
sufficiently large. }

\bigskip

\begin{lem}\label{L:etadist}
There is a constant $C_1>1$ independent of $v$ and $j$ so that
\begin{equation}\label{E:etadist}
C_1^{-1}\le\frac{\eta''_{j+1}(v)}{j+1} \left(\dist(v,\partial
J)\right)^2\le C_1
\end{equation}
for $v\in J\setminus\{v_0\}$.
\end{lem}

\begin{proof}
This follows easily from
\begin{equation}\label{E:eta''}
\eta_{j+1}''(v)=
\begin{cases}
(j+1) \sec^2  v,  &-\frac{\pi}{2}<v<v_0;\\
(j+1) \sec^2  (v+\theta), &v_0<v<\frac{\pi}{2}-\theta.
\end{cases}
\end{equation} 
\end{proof}

\bigskip

\begin{lem}\label{L:etavar}
There exists a constant $C_2\ge 1$ independent of
$v^\sharp,j$ such that

\begin{multline*}
{\max}\{ \eta_{j+1}''(v)\; ;\;|v-v^\sharp|\le
\frac12\dist(v^\sharp,\partial J),
\; v\neq v_0\} \\
\le C_2\;
{\min}\{\eta_{j+1}''(v)\;;\;\;|v-v^\sharp|\le
\frac12\dist(v^\sharp,\partial J), \; v\neq v_0\}
\end{multline*}
for all $v^\sharp\in J$.
\end{lem}

\medskip

\begin{proof}
This follows easily from Lemma \ref{L:etadist}, setting $C_2=9C_1^2$.
\end{proof}

\medskip

\begin{proof}[Proof of Proposition \ref{P:est}, part (i)]
We start by noting that 

\begin{multline*}
|F_{j+1}(-x-i h)|\ge\\
\left| \;
\int\limits_{\mu_{j+1}(x)-h_0^{-1}} ^ {\mu_{j+1}(x)+h_0^{-1}}
e^{iv h}\; e^{vx-\eta_{j+1}(v)} dv \right|-\int\limits_{|v-\mu_{j+1}(x)|>
h_0^{-1}} e^{vx-\eta_{j+1}(v)} dv.
\end{multline*}

\medskip

There is no
significant cancellation in the first term of the right hand  side of the above
inequality; in fact we have

\begin{align*}
\left|\,\int\limits_{\mu_{j+1}(x)-h_0^{-1}}^{\mu_{j+1}(x)+h_0^{-1}} e^{iv h} 
e^{vx-\eta_{j+1}(v)} dv \right|
&\ge \Re\int\limits_{\mu_{j+1}(x)-h_0^{-1}}^{\mu_{j+1}(x)+h_0^{-1}}
e^{i(v-\mu_{j+1}(x)) h}e^{vx-\eta_{j+1}(v)} \,dv\\
&\ge \frac{1}{2}
\int\limits_{\mu_{j+1}(x)-h_0^{-1}}^{\mu_{j+1}(x)+h_0^{-1}}
e^{vx-\eta_{j+1}(v)}\, dv
\end{align*}
since $\cos \alpha>\frac12$ for $|\alpha|\le 1$.

 Let $\epsilon=\epsilon(x,h_0,j)>0$.  We will determine an explicit value of
$\epsilon$ later on, but we will require that
\begin{equation}\label{E:epsdist}
\epsilon(x,h_0,j)\le\frac12\dist(\mu_{j+1}(x),\partial J)
\end{equation}
and
\begin{equation}\label{E:epsw}
\epsilon(x,h_0,j)\le\frac{1}{2 h_0}.
\end{equation}

\medskip 

Thus in particular

\begin{align*}
|F_{j+1}(-x-i h)| &\ge
\frac{1}{2}\;\int\limits_{\mu_{j+1}(x)-h_0^{-1}}^{\mu_{j+1}(x)+h_0^{-1}}
e^{vx-\eta_{j+1}(v)} \,dv-\\
&\qquad\qquad-\int\limits_{|v-\mu_{j+1}(x)|> h_0^{-1}} e^{v
x-\eta_{j+1}(v)} \,dv\\
&\ge\frac{1}{2}
\int\limits_{\mu_{j+1}(x)-\epsilon}^{\mu_{j+1}(x)+\epsilon}{ e^{v x-
\eta_{j+1}(v)} \,dv}-\\
&\qquad\qquad -\int\limits_{|v-\mu_{j+1}(x)|> h_0^{-1}}
{e^{vx-\eta_{j+1}(v)}\, dv}.
\end{align*}

\medskip

The integrand is strictly increasing for $v\le \mu_{j+1}(x)$ and strictly
decreasing for $v\ge \mu_{j+1}(x)$; hence

\begin{eqnarray}\label{E:Flow}
|F_{j+1}(-x-i h)| 
&\ge& \frac{\epsilon}{2} \left( e^{\Theta_{j+1}
(\mu_{j+1}(x)-\epsilon)}+e^{\Theta_{j+1}(\mu_{j+1}(x)
+\epsilon)}\right)-\notag\\
&-&|J| \left( e^{\Theta_{j+1}
(\mu_{j+1}(x)-h_0^{-1})}
+e^{\Theta_{j+1}(\mu_{j+1}(x)+h_0^{-1})}\right)
\end{eqnarray}

\noindent where 
\begin{equation}\label{E:Thetadef}
\Theta_{j+1} (v):= vx-\eta_{j+1}(v).
\end{equation}

\bigskip
Let us work for the moment  with the case that $$v_{0}\notin
\left(\mu_{j+1}(x)-h_0^{-1},
\mu_{j+1}(x)+h_0^{-1}\right).$$

We shall estimate $\Theta_{j+1}$ using the integral form of Taylor's
theorem at $\mu_{j+1}(x)$. We have

{\allowdisplaybreaks
\begin{subequations}\label{E:Taylor}
\begin{align}
\Theta_{j+1}(\mu_{j+1}(x)-\epsilon)&=
\Theta_{j+1}(\mu_{j+1}(x))+\\
&+\int\limits_{\mu_{j+1}(x)-\epsilon}^{\mu_{j+1}(x)}
\Theta_{j+1}''(t)\left(t-(\mu_{j+1}(x)-\epsilon)\right)\,dt\notag\\
\Theta_{j+1}(\mu_{j+1}(x)+\epsilon)&=
\Theta_{j+1}(\mu_{j+1}(x))+\\
&+\int\limits^{\mu_{j+1}(x)+\epsilon}_{\mu_{j+1}(x)}
\Theta_{j+1}''(t)\left(\mu_{j+1}(x)+\epsilon-t\right)\,dt\notag\\
\Theta_{j+1}(\mu_{j+1}(x)-h_0^{-1})&= \Theta_{j+1}(\mu_{j+1}(x))+
\\
&+\int\limits_{\mu_{j+1}(x)-h_0^{-1}}^{\mu_{j+1}(x)}
\Theta_{j+1}''(t)\left(t-(\mu_{j+1}(x)-h_0^{-1})\right)\,dt\notag\\
\Theta_{j+1}(\mu_{j+1}(x)+h_0^{-1})&=
\Theta_{j+1}(\mu_{j+1}(x))+\\
&+\int\limits^{\mu_{j+1}(x)+h_0^{-1}}_{\mu_{j+1}(x)}
\Theta_{j+1}''(t)\left(\mu_{j+1}(x)+h_0^{-1}-t\right)\,dt\notag.
\end{align}
\end{subequations}
}

Using Lemma \ref{L:etavar} to estimate the integrands in (\ref{E:Taylor} a,b) we
 obtain the following estimates:

\begin{subequations}\label{E:Taylorests}
\begin{align}
\Theta_{j+1}(\mu_{j+1}(x)-\epsilon)&\ge 
\Theta_{j+1}(\mu_{j+1}(x))-\frac{C_2
\epsilon^2}{2}\;\eta_{j+1}''(\mu_{j+1}(x))
\tag{\ref{E:Taylorests}a}\\
\Theta_{j+1}(\mu_{j+1}(x)+\epsilon)&\ge 
\Theta_{j+1}(\mu_{j+1}(x))-\frac{C_2
\epsilon^2}{2}\;\eta_{j+1}''(\mu_{j+1}(x))
\tag{\ref{E:Taylorests}b}.
\intertext{ Reducing the intervals of integration in
\eqref{E:Taylor} (c) and (d) to $[\mu_{j+1}(x)-\epsilon,\mu_{j+1}(x)]$ and
$[\mu_{j+1}(x),\mu_{j+1}(x)+\epsilon]$, respectively, and invoking Lemma
\ref{L:etavar} and the assumption \eqref{E:epsw} we see that the integrands are
bounded above by $-\frac{\eta_{j+1}''(\mu_{j+1}(x))}{2 C_2 h_0}
$.  This yields}
\Theta_{j+1}(\mu_{j+1}(x)-h_0^{-1})&\le
\Theta_{j+1}(\mu_{j+1}(x))-\frac{\epsilon}{2 C_2  h_0}\,
\eta_{j+1}''(\mu_{j+1}(x))
\tag{\ref{E:Taylorests}c}\\
\Theta_{j+1}(\mu_{j+1}(x)+h_0^{-1})&\le
\Theta_{j+1}(\mu_{j+1}(x))-\frac{\epsilon}{2 C_2 h_0}\,
\eta_{j+1}''(\mu_{j+1}(x)).
\tag{\ref{E:Taylorests}d}
\end{align}
\end{subequations}
(If $\mu_{j+1}(x)- h_0^{-1}$ or $\mu_{j+1}(x)+ h_0^{-1}$ lands outside of $J$
then the corresponding estimate holds by default.)

\medskip

\medskip

If $v_0\in\left(\mu_{j+1}(x)-h_0^{-1},\mu_{j+1}(x)\right)$ then 
(\ref{E:Taylor}c) must be adjusted by inclusion of the term
\begin{align*}
&\left(\Theta'_{j+1}(v_0+)-\Theta'_{j+1}(v_0-)\right)
\left(v_0-(\mu_{j+1}(x)-h_0^{-1})\right)\\
&\qquad\qquad=\left(-\eta'_{j+1}(v_0+)+\eta'_{j+1}(v_0-)\right)
\left(v_0-(\mu_{j+1}(x)-h_0^{-1})\right)\\
&\qquad\qquad\le 0.
\end{align*} 
Thus (\ref{E:Taylorests}c),  still holds, as well as (\ref{E:Taylorests}b) and
(\ref{E:Taylorests}d), though (\ref{E:Taylorests}a) may fail.

Similar considerations apply to the case where
$v_0$ lies in the interval $\left(\mu_{j+1}(x),\mu_{j+1}(x)+h_0^{-1}\right)$.

Applying all of this to \eqref{E:Flow} (and dropping the term for which we have
no positive lower bound)  we obtain:

\begin{align}\label{E:key}
&|F_{j+1}(-x-i h)|\notag\\
&\qquad \ge e^{\Theta_{j+1}(\mu_{j+1}(x))} 
\;\left( \frac{\epsilon}{2}\,e^{-\frac{C_2
\epsilon^2}{2}\;\eta_{j+1}''(\mu_{j+1}(x))}-2 |J|
\,e^{\frac{-\epsilon}{2 C_2 h_0}\,
\eta_{j+1}''(\mu_{j+1}(x))}\right)\notag\\
&\qquad =
e^{\tilde\eta_{j+1}(x)} %%%%e^{\widetilde{\eta_{j+1}}(x)} 
\;\left( \frac{\epsilon}{2}\,e^{-\frac{C_2
\epsilon^2}{2}\;\eta_{j+1}''(\mu_{j+1}(x))}-2 |J|
\,e^{-\frac{\epsilon}{2 C_2 h_0}\,
\eta_{j+1}''(\mu_{j+1}(x))}\right).
\end{align}

\medskip

We shall choose $\epsilon$ such that
\begin{equation*}\label{E:nointerfere}
2 |J|
\,e^{-\frac{\epsilon}{2 C_2 h_0}\,
\eta_{j+1}''(\mu_{j+1}(x))}\le 
\frac{\epsilon}{4} \,e^{-\frac{C_2 \epsilon^2}{2}\;\eta_{j+1}''(\mu_{j+1}(x))},
\end{equation*}
\noindent
i.e. 
\begin{equation}\label{E:lastepscond}
\log \frac{ 8|J|}{\epsilon}+\frac{C_2 \epsilon^2}{2}\;\eta_{j+1}''(\mu_{j+1}(x))
\le \frac{\epsilon}{2 C_2 h_0}\,
\eta_{j+1}''(\mu_{j+1}(x)).
\end{equation}

Then
\begin{equation}\label{E:punchi}
|F_{j+1}(-x-i h)|\ge
\frac{\epsilon}{4}\,e^{-\frac{C_2
\epsilon^2}{2}\;\eta_{j+1}''(\mu_{j+1}(x))}
 e^{\tilde\eta_{j+1}(x)}. 
%%%%%%e^{\widetilde{\eta_{j+1}}(x)}.
\end{equation}

\medskip
It remains to choose $\epsilon$ satisfying \eqref{E:epsdist}, \eqref{E:epsw}, 
\eqref{E:lastepscond}.  We can choose 
\begin{equation}\label{E:epsdef}
\epsilon:= \frac{C_3}{\sqrt{\eta_{j+1}''(
\mu_{j+1}(x))}}
\end{equation}
with $C_3=C_3(h_0)>0$ independent of $x$ and $j$.

Since $\eta_{j+1}''(
\mu_{j+1}(x))$ is bounded below uniformly in $j$, \eqref{E:epsw} will hold if
$C_3$ is small enough.  

Lemma \ref{L:etadist} shows that \eqref{E:epsdist} is also guaranteed for small
$C_3$.

Condition \eqref{E:lastepscond} now reads
\begin{equation*}\label{E:etabig}
\log \left(\frac{ 8|J|}{C_3}\sqrt{\eta_{j+1}''(\mu_{j+1}(x))}\right)
+\frac{C_2C_3^2}{2}
\le \frac{C_3}{2 C_2 h_0}\,
\sqrt{\eta_{j+1}''(\mu_{j+1}(x))}.
\end{equation*}

This will hold provided that $\eta_{j+1}''(\mu_{j+1}(x))$ exceeds some absolute
constant $M$.

Consulting \eqref{E:eta''} we see that $\eta_{j+1}''(\mu_{j+1}(x))>M$ in the
following cases:
\begin{itemize}
\item for any fixed $j$ provided that $x$ is large enough;
\item for all $x$ provided that $j$ is large enough.
\end{itemize}

Combining \eqref{E:epsdef} with \eqref{E:punchi} we see that in these cases we
have
\eqref{E:noncorner} with $L=\frac{C_3}{4}e^{-C_2C_3^2/2}$.
\end{proof}

\medskip
\begin{proof}[Proof of Proposition \ref{P:est}, part (ii)]
We explain where the proof of part (i) must be modified.

In the current case we have $\mu_{j+1}(x)=v_0$, so
$\Theta_{j+1}$ will not be differentiable at $\mu_{j+1}(x)$, but the one-sided
derivatives $\Theta'_{j+1}(v_0+)$ and $\Theta'_{j+1}(v_0-)$ will exist.

Recalling \eqref{E:Thetadef} we have 
\begin{align*}\label{E:jump}
\Theta'_{j+1}(v_0+)-\Theta'_{j+1}(v_0-)=-(j+1)B,
\end{align*}
where $B:=\eta'_{1}(v_0+)-\eta'_{1}(v_0-)>0$.

Since $\Theta_{j+1}$ has a maximum at $v_0$, we must have
\begin{equation}\label{E:maxcond}
\Theta'_{j+1}(v_0+)\le0\le \Theta'_{j+1}(v_0-).
\end{equation}

Combining \eqref{E:Thetadef} and \eqref{E:maxcond} we find that
\begin{equation}\label{E:jumpbnd}
-(j+1)B \le \Theta'_{j+1}(v_0+)\le 
0 \le \Theta'_{j+1}(v_0-)\le (j+1)B.
\end{equation}

The Taylor expansions \eqref{E:Taylor}  must be modified by inclusion on the
right-hand side of the terms $-\Theta'_{j+1}(v_0-)\epsilon,$
$\Theta'_{j+1}(v_0+)\epsilon,$
$-\Theta'_{j+1}(v_0-)h_{0}^{-1},$ and $\Theta'_{j+1}(v_0+)h_{0}^{-1},$ respectively.

Focusing on the latter two expansions, we see that the new terms are negative.
If we restrict the  integrals to
$[v_0-\min\{h_{0}^{-1},\frac12\dist(v_0,\partial J)\}, v_0]$ and
$[v_0,v_0+\min\{h_{0}^{-1},\frac12\dist(v_0,\partial J)\}]$, respectively, we
obtain   modified versions of (\ref{E:Taylorests} c,d) taking the following
form:

\begin{subequations}\label{E:modcd}
\begin{align}
\Theta_{j+1}(\mu_{j+1}(x)-h_0^{-1})&\le
\Theta_{j+1}(\mu_{j+1}(x))-(j+1)C_4 
\tag{\ref{E:modcd}c}\\
\Theta_{j+1}(\mu_{j+1}(x)+h_0^{-1})&\le
\Theta_{j+1}(\mu_{j+1}(x))-(j+1)C_4,
\tag{\ref{E:modcd}d}
\end{align}
\end{subequations}

\medskip

\noindent where $C_4$ depends on $h_0$ but not on $j$ or $x$.

To obtain modified versions of  (\ref{E:Taylorests}a,b) we apply the inequalities 
\eqref{E:jumpbnd} to  the new first derivative terms to obtain the following:

\begin{align*}
\Theta_{j+1}(\mu_{j+1}(x)-\epsilon)&\ge 
\Theta_{j+1}(\mu_{j+1}(x))-(j+1)B\epsilon-(j+1)C_5\epsilon^2
\tag{\ref{E:modcd}a}\\
\Theta_{j+1}(\mu_{j+1}(x)+\epsilon)&\ge 
\Theta_{j+1}(\mu_{j+1}(x))-(j+1)B\epsilon-(j+1)C_5\epsilon^2.
\tag{\ref{E:modcd}b}
\end{align*}

The inequality \eqref{E:lastepscond} is now modified to read
\begin{equation}\label{E:modepscond}
\log\frac{4|J|}{\epsilon} +(j+1)B\epsilon+(j+1)C_5\epsilon^2
\le (j+1)C_4.
\end{equation}

If we now set $\epsilon=\frac{C_6}{j+1}$ ($C_6>0$ depends on $h_0$ but not
on
$x$ or $j$) we find that \eqref{E:epsdist} and \eqref{E:epsw} will hold if $C_6$
is chosen small enough, while \eqref{E:modepscond} holds for all large enough
$j$.

Under these conditions we obtain the following modified form of
\eqref{E:punchi}:

\begin{align*}\label{E:punchii}
|F_{j+1}(-x-i h)|&\ge
\frac{\epsilon}{2}\,e^{-(j+1)(B\epsilon+C_5 \epsilon^2)}
e^{\tilde\eta_{j+1}(x)}\\ %%%%e^{\widetilde{\eta_{j+1}}(x)} \\
&\ge\frac{K}{j+1}e^{\tilde\eta_{j+1}(x)}, %%%%e^{\widetilde{\eta_{j+1}}(x)}
\end{align*}
with $K=\frac{C_6}{2}e^{-B C_6-C_5 C_6^2} $.
\end{proof}

\begin{coro}\label{C:finstrip}
The union of the zero sets of the $F_{j+1}$ is finite (counting multiplicity) in
any strip
$-h_0\le\Im\xi\le 0$.
\end{coro}

\begin{proof}
By part (ii) of Proposition \ref{P:est}, the $F_{j+1}$ are zero-free when
$j\ge j_0(h_0,|J|)$.  But by the remark in the statement of Proposition
\ref{P:est}, the zero set of $F_{j+1}$ in our strip is finite for $0\le
j<j_0(h_0,|J|)$.
\end{proof}

From Proposition \ref{P:est} we see that the union of the zero sets of the $F_j$
contains finitely many points in each strip $-h_0\le \Im\xi\le 0$.  Thus in
particular, for all but a discrete set of $h>0$ the $F_j$ are all non-vanishing on
the line $\Im \xi=-h$.  We assume for the rest of this section that $h$ has been
chosen with the property.

Returning to the integrals \eqref{E:keyint}, we see from Proposition  \ref{P:est},
\eqref{E:Legt} and \eqref{E:Imtbnd} that the integrand decays exponentially on
the strip $-h\le \Im\xi\le 0$.  Thus we may apply the residue theorem on
this strip as indicated in \S\ref{S:Bker}. We still need upper bounds for the
shifted integrals 
\begin{equation*}\label{E:shifty}
\int\limits^{\infty}_{-\infty} \frac{e^{i(t_2-\overline{\tau}_2)\, 
\frac{(-x-ih)}{2}}\;dx}{F_{j+1}(-x-i h)}.
\end{equation*}

Setting $b:=\Im(t_2-\overline{\tau}_2)\in (-\pi,\pi-2\theta)$,  we have
\begin{equation}\label{E:shiftrev}
\left| \int\limits^{\infty}_{-\infty} \frac{e^{i(t_2-\overline{\tau}_2)\, 
\frac{(-x-ih)}{2}}\;dx}{F_{j+1}(-x-i h)}\right|
\le e^{h\Re(t_2-\bar\tau_2)/2}
\int\limits_{-\infty}^\infty 
\frac{e^{ x \frac{b}{2}}}{\left|F_{j+1}(-x-i h)\right|} dx.
\end{equation}

From Proposition  \ref{P:est} together with \eqref{E:Legt}, \eqref{E:muspec},
\eqref{E:eta''} we have
\begin{multline*}
\frac{e^{ x \frac{b}{2}}}{\left|F_{j+1}(-x-i h)\right|}
\le\\
\begin{cases}
\frac{1}{L}
\sqrt{(1+j)\left(1+\frac{x^2}{(1+j)^2}\right)}e^{x\frac{b}{2}-
\tilde\eta_{j+1}(x)}%%%%
\\
\hspace{1.0in}\text{for }x\notin[(j+1)\tan v_0,(j+1)\tan (v_0+\theta)];\\
\frac{j+1}{K}e^{x\frac{b}{2}-
\tilde\eta_{j+1}(x)}%%%%
\\
\hspace{1.0in}\text{for } x\in[(j+1)\tan v_0,(j+1)\tan (v_0+\theta)]\\
\end{cases}
\end{multline*}
for $j$ large.

Set
\begin{equation}\label{E:Hdef}
H:=-\inf\{\widetilde\eta_1(x); x\in\RR\}+1.
\end{equation} 

\begin{lem}\label{L:epstrick}
For all $\delta\in(0,1)$ there is $M=M(\delta)$ independent of $j, x$ so that 
\begin{equation*}
\frac{e^{ x \frac{b}{2}}}{\left|F_{j+1}(-x-i h)\right|}
\le M 
e^{x\frac{b}{2}+(\delta-1)\tilde\eta_{j+1}(x)+\delta (j+1)H}.%%%%
\end{equation*}
\end{lem}

\begin{proof}
Setting $\alpha=\frac{x}{j+1}$ and recalling that 
$\tilde \eta_{j+1}(\zeta)=(j+1) \tilde\eta_1\left( \tfrac{\zeta}{j+1}\right)$
 we must choose $M$ so that 
\begin{align*}
\frac{j+1}{K}\le
Me^{(j+1)\delta( \tilde\eta_1(\alpha)+H)} &\text{ for }\alpha\in[\tan v_0,\tan
(v_0+\theta)],%%%%
\\
\frac{\sqrt{j+1}}{L}\sqrt{1+\alpha^2}
\le Me^{(j+1)\delta(\tilde\eta_1(\alpha)+H)}&\text{ for }\alpha\notin[\tan
v_0,\tan (v_0+\theta)].%%%%
\end{align*} 
  It will suffice to choose $M$ so
that
\begin{equation}\label{E:M1}
\left(\frac{1}{L^2}+\frac{1}{K}\right)(j+1)\le
Me^{(j+1)\delta(\tilde\eta_1(\alpha)+H)}%%%%
\end{equation}
and
\begin{equation}\label{E:M2}
1+\alpha^2\le Me^{(j+1)\delta(\tilde\eta_1(\alpha)+H)}%%%%
\end{equation}
hold for all $\alpha$.  That \eqref{E:M1} is possible follows from the fact that
the right-hand side exceeds $Me^{(j+1)\delta}$.  Similarly, \eqref{E:M2} is
possible since the right-hand side exceeds
$Me^{\delta(\tilde\eta_1(\alpha)+H)}$, which grows exponentially with
$\alpha$.
\end{proof}

\begin{lem}\label{L:qscon} Let $X$ be a compact subset of $J$.  Then there  are
$Q>0$, $\gamma=\gamma(X)>0$  independent of $j$ so that 
\begin{align*}
\widetilde\eta'_{j+1}(x-(j+1)Q)\le
\widetilde\eta'_{j+1}(x)-\gamma\\
\widetilde\eta'_{j+1}(x+(j+1)Q)\ge
\widetilde\eta'_{j+1}(x)+\gamma
\end{align*}
when $\widetilde\eta'_{j+1}(x)\in X$.
\end{lem}

\begin{proof}
Again setting $\alpha=\frac{x}{j+1}$ we are reduced to the case $j=0$.  Picking
$Q>\tan(v_0+\theta)-\tan v_0$, our claim follows easily from the fact that
$\widetilde \eta'_1:\mathbb R\to J$ is a continuous non-decreasing surjective
function which is strictly increasing off of the interval $[\tan v_0,\tan(v_0+\theta)]$.
\end{proof}

Using the duality theorem $\tilde{\tilde \eta}=\eta$ (see [H\"o, Thm. 1.3.3]) we
have
\begin{align*}
&\sup\left\{x\frac{b}{2}-(1-\delta)\tilde\eta_{j+1}(x)\,; x\in\RR\right\}\\
&\qquad = (1-\delta)\sup\left\{x\frac{b}{2(1-\delta)}-\tilde\eta_{j+1}(x)\,;
x\in\RR\right\}\\
&\qquad =(1-\delta)
\widetilde{\widetilde{\eta_{j+1}}}\left(\frac{b}{2(1-\delta)}\right)\\ 
&\qquad=(1-\delta) \eta_{j+1}\left(\frac{b}{2(1-\delta)}\right).
\end{align*}

Since $ x \frac{b}{2}-\widetilde{\eta}_{j+1}(x) \to -\infty$ as $|x|\to +\infty$ and $\widetilde{\eta}_{j+1}$ is $C^1 (\mathbb{R})$ we can find an  $x_0$ so that
$x_0\frac{b}{2}-(1-\delta)\widetilde\eta_{j+1}(x_0)=(1-\delta)
\eta_{j+1}\left(\frac{b}{2(1-\delta)}\right)$. But then  $\frac{d}{dx}
\left(x\frac{b}{2}-(1-\delta)\widetilde\eta_{j+1}(x)\right)$ vanishes at
$x=x_0$, so we can deduce from Lemma \ref{L:qscon} and the integral form of 
Taylor's theorem  that
\begin{multline*}
x\frac{b}{2}-(1-\delta)\widetilde\eta_{j+1}(x)
\le \\
\begin{cases}
(1-\delta)
\eta_{j+1}\left(\frac{b}{2(1-\delta)}\right)+(1-\delta)\gamma(x-(x_0-(j+1)Q))\\
\hspace{1in}\text{ for } x\le x_0-(j+1)Q;\\
(1-\delta) \eta_{j+1}\left(\frac{b}{2(1-\delta)}\right)
-(1-\delta)\gamma(x-(x_0+(j+1)Q))\\
 \hspace{1in}\text{ for } x\ge x_0+(j+1)Q
\end{cases}
\end{multline*}
provided that $\widetilde\eta'_{j+1}(x_0)=\frac{b}{2(1-\delta)}\in X$.

Combining this with Lemma \ref{L:epstrick} we have
\begin{align}\label{E:Legdual}
&\int\limits_{-\infty}^\infty 
\frac{e^{ x \frac{b}{2}}}{\left|F_{j+1}(-x-i h)\right|} dx\\
&\qquad \le
Me^{(1-\delta)\eta_{j+1}\left(\frac{b}{2(1-\delta)}\right)+\delta(j+1)H}
\int\limits_{-\infty}^{x_0-(j+1)Q}
e^{(1-\delta)\gamma(x-(x_0-(j+1)Q))}\,dx +\notag\\
&\qquad +Me^{(1-\delta)\eta_{j+1}\left(\frac{b}{2(1-\delta)}\right)+\delta(j+1)H}
\int\limits_{x_0-(j+1)Q}^{x_0+(j+1)Q}\,dx+\notag\\
&\qquad +Me^{(1-\delta)\eta_{j+1}\left(\frac{b}{2(1-\delta)}\right)+\delta(j+1)H}
\int\limits_{x_0+(j+1)Q}^{\infty}
e^{-(1-\delta)\gamma(x-(x_0+(j+1)Q))}\,dx\notag\\
&\qquad  = M\left(2(j+1)Q+\frac{2}{(1-\delta)\gamma}\right)
e^{(1-\delta)\eta_{j+1}\left(\frac{b}{2(1-\delta)}\right)+\delta(j+1)H}.\notag
\end{align}

If $\frac{b}{2}$ is restricted to a compact subset $X$ of $J$ then 
$\frac{b}{2(1-\delta)}$ is restricted to a slightly larger compact subset $X'$
provided that $\delta\le\delta_0(X)$.

Combining \eqref{E:Legdual} with \eqref{E:shiftrev}, we have proved the following.

\begin{proposition}\label{E:upb}
 For $X$ a compact subset of $J$ there are $R=R(X)>0$, $\delta_0(X)>0$ so that
\begin{equation*}
\left| \int\limits_{{\mathbb R}-ih} \frac{e^{i(t_2-\overline{\tau}_2)\,
\frac{\xi}{2}}\;d\xi} {F_{j+1}(\xi)}\right|
\le
R\cdot(j+1)e^{(1-\delta)\eta_{j+1}\left(\frac{b}{2(1-\delta)}\right)+\delta(j+1)H}
\end{equation*}
when $\frac{b}{2}=\frac12\Im(t_2-\overline{\tau}_2)\in X$ and
$0<\delta<\delta_0(X)$.
\end{proposition}

With harder work one can obtain the following sharper estimates in Proposition
\ref{E:upb}: For $X$ a compact subset of $J$ there exists  $R'=R'(X)>0$ so that

\begin{equation*}
\left| \int\limits_{{\mathbb R}-ih} \frac{e^{i(t_2-\overline{\tau}_2)\,
\frac{\xi}{2}}\;d\xi} {F_{j+1}(\xi)}\right|
\le
R'\cdot(j+1)^2 \; e^{\eta_{j}(\tfrac{b}{2})}
\end{equation*}
when $\frac{b}{2}=\frac12\Im(t_2-\overline{\tau}_2)\in X$.

\section{The asymptotic 
formula and regularity for the Bergman kernel in $\hat\Omega'$ }
\label{S:reg}

\medskip
\vspace{1cm}

Let  $\eta(\frac{b}{2}):=\eta_{1}(\frac{b}{2})=-\log\psi_{r,\theta}(\frac{b}{2}).$

\medskip

Suppose that we choose $t=(t_1,t_2),\;(\tau_1, \tau_2) \in D'$ such that

\bea
|t_1|^2 &\le& e^{-\eta(u)}\\
|\tau_1|^2 &<& e^{-\eta(v)}.\\
\eea

\noindent for $u=\Im t_2\in\overline J,  v=\Im\tau_2\in J$.

Set $b=u+v$ as before.
By the convexity of $\eta$ we have:

$$
|t_1|\; |\tau_1| \;e^{\eta(\frac{b}{2})}\le
|t_1|\;|\tau_1|\;e^{\frac{\eta(u)+\eta(v)}{2} }< 1.
$$

Choose $\delta>0$ small so that
\begin{equation}\label{E:delcheat}
|t_1||\tau_1|e^{(1-\delta)\eta\left(\frac{b}{2(1-\delta)}\right)+\delta H}<1,
\end{equation}
where $H$ is defined in \eqref{E:Hdef}.

\medskip

Consider the series:

\begin{equation}\label{E:error}
\sum_{j\ge 0} \frac{j+1}{4 \pi^2} (t_1 \;\overline{\tau_1})^{j}
\int\limits_{\mathbb R}
\frac{e^{i(t_2-\overline{\tau}_2) \frac{(-x-i h)}{2}}\;\; dx}{ F_{j+1}(-x-i h)}\;.
\end{equation}

\medskip

Setting $t_2-\overline{\tau_2}:=a+i b$ and invoking Proposition \ref{E:upb},  the 
absolute value of the above series can be majorized by

$$
R\; e^{a \frac{h}{2}}\;\sum_{j\ge 0} (j+1)^{2} |t_1
\;\overline{\tau_1}|^{j}  
 \; e^{(j+1)\left((1-\delta)\eta\left(\frac{b}{2(1-\delta)}\right)+\delta H\right)}.
$$

\medskip

Applying the root test  and taking into account \eqref{E:delcheat} we have

\begin{align*}
\limsup&\sqrt[j]{(j+1)^{2} |t_1
\;\overline{\tau_1}|^{j}  
 \; e^{(j+1)\left((1-\delta)\eta\left(\frac{b}{2(1-\delta)}\right)+\delta H\right)}}\\
&=|t_1 \overline{\tau}_1| \;\;
e^{\left((1-\delta)\eta\left(\frac{b}{2(1-\delta)}\right)+\delta H\right)}\\
&< 1.
\end{align*}

\medskip

Thus  we have proved the following:
\begin{prop}\label{P:sum}
The series \eqref{E:error} is $O\left(e^{h\Re t_2/2}\right)$ as $\Im
t_2\to-\infty$, uniformly as $\tau$ ranges over any compact subset of $D'$.
\end{prop}

\medskip

{\bf \centerline  {Positive regularity results}} 

\medskip

Choose $h>0$ so that the strip $-h\le\Im\xi\le 0$ is zero-free for all 
$F_{j+1}(\xi)$.

Then for appropriate $t,\tau$ we have:

$$
K_{D'} (t, \tau)=\sum_{j\ge 0} \frac{j+1}{4 \pi^2} (t_1 \overline{\tau_1})^{j} 
\int\limits_{\mathbb R} \frac{ e^{i(t_2-\overline{\tau_2}) 
\frac{(-x-i h)}{2}}}{F_{j+1}(-x-i h)} dx=
O\left(e^{h\Re t_2/2}\right).
$$

\medskip

Let us fix $\tau \in  D'$. 
A careful inspection of our earlier work shows that we may differentiate the above
formula with respect to
$t_1, t_2$  to obtain:

$$
D^{l}_{t_1}K_{D'}(t,\tau)=\sum_{j\ge l} \tfrac{ (j+1) j (j-1)\cdots (j-l+1)}{4 
\pi^2} \,\overline{\tau_{1}}^{l}\,
(t_1 \overline{\tau_1})^{j-l} \int\limits_{\mathbb R} \frac{
e^{i(t_2-\overline{\tau_2})
\frac{(-x-i h)}{2}}}{ F_{j+1}(-x-ih)} dx
$$

\noindent and 

$$
D^{l}_{t_2} K_{D'} (t,\tau)=\sum_{j\ge 0} 
\frac{j+1}{4 \pi^2} (t_1 \overline{\tau_1})^{j} \int\limits_{\mathbb R} \frac{
e^{i(t_2-\overline{\tau_2}) \frac{(-x-i h)}{2}}
 \left(\frac{i(-x-i h)}{2}\right)^{l}}{ F_{j+1}(-x-i h)}\,
dx
$$

\noindent with corresponding formulae for mixed partials.

\medskip
Arguing as in Proposition \ref{P:sum} we find that each partial derivative
satisfies
\begin{equation}\label{E:mpest}
D^{l_1}_{t_1}D^{l_2}_{t_2} K_{D'} (t,\tau)=O\left(e^{h\Re t_2/2}\right)\text{ as
}\Im t_2\to -\infty
\end{equation}

\medskip 

Let us fix a point $\tilde\omega \in \hat\Omega'$. 
Using the transformation formula \eqref{E:Btran} for the Bergman kernel 
and differentiating with respect to $\tilde w_1, \tilde w_2$  we
obtain:

\begin{align*}
\frac{\partial}{\partial \tilde w_1} K_{\hat\Omega'}(\tilde w, \tilde{\omega})
&= \frac{1}{2 ( \tilde w_2 \overline{\tilde\omega_2})^{\frac{3}{2}}}
\frac{\partial}{\partial t_1} K_{D'}( \Phi(\tilde w), \Phi(\tilde\omega)) 
\frac{\partial
t_1}{\partial \tilde w_1}\\
\frac{\partial}{\partial \tilde w_2} K_{\hat\Omega'}(\tilde w, \tilde{\omega})
&=
-\frac{3}{4 (\tilde w_2)^{\frac{5}{2}}
 (\overline{\tilde{\omega_2})}^{\frac{3}{2}}} K_{D'}( \Phi(\tilde w),
\Phi(\tilde\omega))+\\
&+\frac{1}{2 ( \tilde w_2 \overline{\tilde\omega_2})^{\frac{3}{2}}}\left(
\frac{\partial t_1}{\partial \tilde w_2}
\frac{\partial }{\partial t_1} K_{D'}( \Phi(\tilde w), \Phi(\tilde\omega))\right)+\\
&+\frac{1}{2 ( \tilde w_2 \overline{\tilde\omega_2})^{\frac{3}{2}}}\left(
\frac{\partial t_2}{\partial \tilde w_2}
\frac{\partial }{\partial t_2} K_{D'}( \Phi(\tilde w), \Phi(\tilde\omega))\right).
\end{align*}

%\frac{\partial t_2}{\partial \tilde w_2} \frac{\partial}{\partial t_2} K_{D'}
%( \Phi(\tilde w), \Phi(\tilde\omega))).
%\end{align*}

\medskip 
But 

\bea
\frac{\partial t_1}{\partial \tilde w_1}&=& \frac{1}{ \sqrt{ 2\, \tilde w_2}},\\
\frac{\partial t_1}{\partial \tilde w_2}&=& \frac{-\tilde w_1}{2 \sqrt{2 \,\tilde w_2}} \frac{1}{\tilde w_2},\\
\frac{\partial t_2}{\partial \tilde w_2}&=& \frac{1}{\tilde w_2}.\\
\eea

Combining the transformation laws with \eqref{E:mpest} 
and the fact that $|\tilde w_1|
\le |\sqrt{\tilde w_2}|$ in
$\overline{\hat\Omega'}$  we see that 

$$
|\nabla_{\tilde w} K_{\hat\Omega'}(\tilde w, \tilde\omega)| \le 
O\left({\tilde w_2 }^{\frac{h-5}{2}}\right)
$$

\noindent  and, more generally,

\begin{equation}\label{G: NablaDrvt}
|\nabla^k_{\tilde w} K_{\hat\Omega'}(\tilde w, \tilde\omega)| \le 
O\left({\tilde w_2 }^{\frac{h-3-2k}{2}}\right)
\end{equation}

\noindent for $\tilde w$ near $0$.

\medskip
Let $\tilde \omega$ be a point inside $\hat \Omega'$. We are interested in the 
regularity of $K_{\hat{\Omega'}}(\cdot,\;\tilde \omega)$ near the complex 
tangent point $(0,0)$. Since $\hat\Omega'$ is Lipschitz and $K(\cdot, \tilde \omega)
$ is harmonic we shall estimate the Sobolev $L^p_{k+\epsilon}$ norm and Besov
$B^{\infty}_{k+\epsilon}$ norm 
of $K(\cdot,\tilde \omega)$ in a neighborhood of the complex tangent point 
using the following theorems by Jerison-Kenig  [JK]: 

\begin{theo} \label{T:JK} A) Let $\Omega$ be a bounded Lipschitz domain in $\RR^{n}$. Let 
$\delta(x)$ be the distance of $x$ from the boundary of $\Omega$. Define
$\nabla^{k} u$ as the vector of all $k^{th}$ order derivatives of a function 
$u$. Suppose that $u$ is a harmonic function in $\Omega$. Let $0\le \epsilon
\le 1$, let $k$ be a nonnegative integer, and let $1<p<\infty$. Then the 
following are equivalent:

\bea
&i)&  \; u \;\;{\mbox{belongs  to}}\;\; L^{p}_{k+\epsilon}(\Omega),\\
\newline
&ii)& \;\; \delta^{1-\epsilon}\,| \nabla^{k+1} u|+ | \nabla^{k} u|+|u|\;\;
{\mbox{belongs  to}}\;\; L^{p}(\Omega).\\
\eea

\medskip
B) Suppose that $u$ is a harmonic function in $\Omega$. Let $0< \epsilon
<1$, let $k$ be a nonnegative integer, and let $1\le p\le \infty$. Then the 
following are equivalent:

\bea
&i)&  \; u \;\;{\mbox{belongs  to}}\;\; B^{p}_{k+\epsilon}(\Omega),\\
\newline
&ii)& \;\; \delta^{1-\epsilon}\,| \nabla^{k+1} u|+ | \nabla^{k} u|+|u|\;\;
{\mbox{belongs  to}}\;\; L^{p}(\Omega).\\
\eea
\end{theo}

\medskip 

To apply these results to $\hat\Omega'$ we use the following facts (valid for
$1<p<\infty$):

\medskip

\begin{enumerate}
\item[(a)] $\tilde w_1^{k} \tilde w_2^d \in L^{p}$ (in a neighborhood of the 
complex tangent point (0,0)) if and only if $\Re d+\tfrac{k}{2}+\tfrac{3}{p}>0$, for $k\in\mathbb{N}, \;
d\in \mathbb{C}$;
\item[(b)] $\delta^{1-\sigma} \tilde w_1^{k} \tilde w_2^{d} \in L^p$ (in a neighborhood
of the point $(0,0)$)  if and only if $\Re d+\tfrac{k}{2}+\tfrac{3}{p}+1 > \sigma,$ for all $\sigma\;
\text{with}  \; 0\le \sigma\le 1, \; k\in \mathbb{N}, \; d\in \mathbb{C}$; here $\delta=\delta(\tilde w,
\partial{\hat \Omega'})$ is the distance to the boundary of $\hat \Omega'$. 

\end{enumerate}

\medskip

Let $1<p<\infty$. Suppose that $\frac{h-3}{2}+\frac{3}{p}>0$.  Using (b) above in combination with
Theorem \ref{T:JK}(A) and the inequality \eqref{G: NablaDrvt} we have the following.

\begin{prop}\label{P:posreg}
Let $1<p<\infty$.  If $\frac{h-3}{2}+\frac{3}{p}>0$ and $\tilde\omega\in
\hat\Omega'$ then
$K_{\hat\Omega'}(\cdot,\tilde\omega)\in L^p_s$ in a neighborhood of the complex tangent point $(0,0)$
for 
$0\le s<\frac{h-3}{2}+\frac{3}{p}$.
\end{prop}

\medskip

Note that for $p=2$ we have $K_{\hat\Omega'}(\cdot,\tilde\omega)\in L^2_s$ for $0\le
s<\frac{h}{2}$.

 When $r=1,\; \theta\sim 0$ and positive  then (by \eqref{E:zs}) $h$ can be
chosen to  be any positive number smaller than $3+\frac{4
\theta}{\pi-\theta}$.   

 When $r=1,\; \theta=\pi-\epsilon,\;
\epsilon\sim 0$ and positive, $h$ can be chosen to be very large since the very first zero has imaginary
part $-\frac{4
\pi}{\epsilon}+1$. In this case, $K_{\hat\Omega'}$ will be very regular.  

\medskip
The map $\Psi: B_1 \cap B_2 \to \hat\Omega'$ defined in \eqref{E:Psidef} is a 
diffeomorphism in a neighborhood of the complex tangent point $q$. Using the local 
diffeomorphism
invariance of Sobolev spaces (see for example [T], Chap. XI,\S 2) we can conclude
that 
$K_{B_1 \cap B_2}\in L^{p}_{s}$ in a neighborhood of $q$ for  
$0\le s<\frac{h-3}{2}+\frac{3}{p}$.

\bigskip

{\bf \centerline{ Negative regularity results}}

\medskip

We begin by investigating the regularity of terms $\tilde w_1^j \tilde w_2^a$,
assuming $j\in\NN\cup\{0\}$, $a\notin \NN\cup\{0\}$.  We claim
that
$\tilde w_1^j \tilde w_2^a \notin L^p_s$ for $s\ge\max\{0,\Re a +
\frac{j}{2}+\frac{3}{p}\}$.  To see this, write $s=l+\sigma$ with 
$l\in\NN\cup\{0\}$, $0\le\sigma\le1$.  Then 
$$
\delta^{1-\sigma} (\tfrac{\partial}{\partial \tilde w_2})^{l+1} \left 
( \tilde w_1^{j} \; \tilde w_2^{a}\right) \notin L^p
$$
showing that $\tilde w_1^j \tilde w_2^a$ indeed fails to lie in $L^p_s$.

In particular, $\tilde w_1^{j} \; \tilde{w_2}^{\frac{i
\xi_{j,k}}{2}-\frac{j+3}{2}}\notin L^p_s$ for 
$s\ge\max\{0,\frac{-\Im \xi_{j,k}-3}{2}+\frac{3}{p}\}$, $\frac{i
\xi_{j,k}}{2}-\frac{j+3}{2}\notin \NN\cup\{0\}$.  When $p=2$, our condition on
$s$ simplifies to $s\ge\max\{0,\frac{-\Im \xi_{j,k}}{2}\}$.

\medskip

Consider the special case $r=1, \theta\sim 0$ \, and positive. (The root pattern in
this case was discussed at the end of \S \ref{S:Bker}.) Suppose we choose $h$
such that the  first zero of $F_1$ lies in the strip $-h<\Im \xi< 0$, but no
other residues lie in the closed strip $-h\le \Im \xi\le 0$.
Then the asymptotic formula for the Bergman kernel $K_{\hat \Omega'}$
 will read

$$
K_{\hat \Omega'}(\tilde w, \tilde\omega)= 
C \tilde w_2^{\frac{i \xi_{0,1}}{2}-\frac{3}{2}}
+O\left(\;\tilde w_2^{\frac{h-3}{2}}\;\right),
$$
with $C\ne0$ for most $\tilde\omega$.

\medskip

Taking into account that $ \frac{i\; \xi_{0,1}}{2}= \frac{3}{2}+\frac{2\theta}{\pi-\theta}$ we see that 
\[\tilde{w_2}^{\frac{i \xi_{0,1}}{2}-\frac{3}{2}}\notin
L^{2}_{s}\text{ for }s\ge \frac{3}{2}+\frac{2\theta}{\pi-\theta}.\]
Arguing as above the error term is in $L^2_s$ for $s<\frac{h}{2}$.  Thus 
$K_{\hat
\Omega'}(\tilde w, \tilde\omega)\notin L^{2}_{s}\text{ for }s\ge \frac{3}{2}+
\frac{2\theta}{\pi-\theta}$.

It follows that for every $s>\frac{3}{2}$ we can choose $\hat\Omega'$ and
$\tilde\omega\in \hat\Omega'$ so that $K_{\hat\Omega'}(\cdot,\tilde\omega)\not\in
L^{2}_{s}$.

\medskip

More generally, we have the following for general $r,\theta$.

\begin{prop}\label{P:posreg2}
If the strip $-h_*\le\Im\xi\le 0$ contains 
\begin{enumerate}
\item[a)] a zero $\xi_{j,k}$ of $F_{j+1}$ with $\frac{i
\xi_{j,k}}{2}-\frac{j+3}{2}\notin \NN\cup\{0\}$ 
\end{enumerate}
 or
\begin{enumerate}
\item[b)] a multiple zero of some $F_{j+1}$
\end{enumerate}
and $\frac{h_*-3}{2}+\frac{3}{p}\ge0$
then for most $\tilde\omega\in\hat\Omega'$ we have 
$K_{\hat\Omega'}\left(\cdot,\tilde\omega\right)\notin 
L^p_{\frac{h_*-3}{2}+\frac{3}{p}}$.
\end{prop}

In proving this result, we work on a strip $-h\le\Im \xi\le 0$ with $h$ a
little bit larger than $h_*$.

Applying the same reasoning with $p=\infty$ and applying part (B) rather than
part (A) of Theorem
\ref{T:JK} we find that
$K_{\hat\Omega'}\left(\cdot,\tilde\omega\right)\notin
B^\infty_{s}$ for $s>\frac{h_*-3}{2}$. Recall
that for
$0< \epsilon <1$ the function space $B^{\infty}_{\epsilon}$ coincides  with the
usual H\"older class of order $\epsilon$.

Returning to the special case $r=1, \theta\sim 0$ \, and positive, we find that
for every positive $\epsilon$, we can choose $\hat\Omega'$ and
$\tilde\omega\in \hat\Omega'$ so that $K_{\hat\Omega'}(\cdot,\tilde\omega)$
fails to be H\"older of order $\epsilon$.

As before,  a change of variable argument allows us to transfer all of these
conclusions to the behavior of $K_\Omega(\cdot,\zeta)$ near a complex tangent point
$q\in\partial\Omega$.

\medskip
%\section{Acknowledgments}

%\medskip
%First author supported in part by the National Science
%Foundation under Grant No. DMS-0072237. This work was completed 
%while the second author was visiting the University of Toronto. 
%She would like to thank the Department of Mathematics for its 
%hospitality and support. 

\end{document}